\documentclass{gtmon_a}
\pdfoutput=1

\proceedingstitle{Exotic homology manifolds (Oberwolfach 2003)}
\conferencestart{29 June 2003}
\conferenceend{5 July 2003}
\conferencename{Workshop on Exotic Homology Manifolds}
\conferencelocation{Oberwolfach Mathematics Institute, Oberwolfach,
Germany}

\editor{Frank Quinn}
\givenname{Frank}
\surname{Quinn}

\editor{Andrew Ranicki}
\givenname{Andrew}
\surname{Ranicki}

\title{Problems on homology manifolds}
\author{Frank Quinn}
\givenname{Frank}
\surname{Quinn}
\address{Department of Mathematics\\Virginia Tech\\Blacksburg\\Virginia
24061-0123\\USA}
\email{quinn@math.vt.edu}
\urladdr{}

\volumenumber{9}
\issuenumber{}
\publicationyear{2006}
\papernumber{6}
\startpage{87}
\endpage{103}

\doi{}
\MR{}
\Zbl{}


\received{17 October 2003} 
\revised{}
\accepted{17 October 2003}
\published{22 April 2006}
\publishedonline{22 April 2006}
\proposed{}
\seconded{}
\corresponding{}
\version{}

\arxivreference{math.GT/0310277}

\makeatletter
\def\cnewtheorem#1[#2]#3{\newtheorem{#1}{#3}[section]
\expandafter\let\csname c@#1\endcsname\c@thm}

\cnewtheorem{lem}[thm]{Lemma}
\cnewtheorem{proposition}[thm]{Proposition}
\cnewtheorem{problem}[thm]{Problem}
\theoremstyle{definition}
\cnewtheorem{defn}[thm]{Definition}

\makeatother  

\renewcommand{\Z}{\ensuremath{\mathbb{Z}}}

\begin{document}

\begin{abstract}
A compilation of questions for the proceedings of the Workshop on Exotic
Homology Manifolds, Oberwolfach June 29 -- July 5, 2003
\end{abstract}

\maketitle

In these notes ``homology manifold'' means ENR (Euclidean neighborhood
retract) $\Z$--coefficient homology manifold, unless otherwise specified,
and ``exotic'' means not a manifold factor (ie local or ``Quinn''
index $\neq1$.  We use the multiplicative version of the local index,
taking values in $1+8\Z$).  In the last decade exotic homology manifolds
have been shown to exist and quite a bit of structure theory has been
developed.  However they have not yet appeared in other areas of
mathematics.  The first groups of questions suggest ways this might
happen.  Later questions are more internal to the subject.

\fullref{S:aspherical}, \fullref{S:limits}, and \fullref{S:compact} concern
possible ``natural'' appearances of
homology manifolds: as aspherical geometric objects; as
Gromov--Hausdorff limits; and as boundaries of compactifications.
\fullref{S:actions} discusses group actions, where the use of homology
manifold fixed sets may give  simpler classification results.
\fullref{S:nonENR} and \fullref{S:approx} consider possible  generalizations to  non-ANR
and ``approximate'' homology manifolds. \fullref{S:special} concerns spaces
with special metric structures. \fullref{S:low} describes still-open low
dimensional cases of the current theory. \fullref{S:DDP}
collects problems related to homeomorphisms and the ``disjoint
disk property'' for exotic homology manifolds.

\section{Aspherical homology manifolds}\label{S:aspherical}
Geometric structures on aspherical spaces seem to be rigid.  The
``Borel conjecture'' is that closed aspherical manifolds are determined
up to homeomorphism by their fundamental groups, and this has been
verified in many cases, see Farrell \cite{F} for a survey.  More
generally it is expected that aspherical homology manifolds should be
determined up to s-cobordism by their fundamental groups, so in
particular the fundamental group should determine the local index.  So
far, however, there are no exotic examples.

\begin{problem}\label{1.1} Is there a closed aspherical homology manifold with local index $\neq1$?\end{problem}
  If
so then exotic homology manifolds would be
required for a full analysis of the aspherical question.
See \fullref{S:approx} for a ``approximate'' version of the problem.

\section{Gromov--Hausdorff Limits}\label{S:limits}
Differential geometers have investigated limits  of Riemannian
manifolds with various curvature and other constraints.  As constraints are weakened one sees
\begin{enumerate}
\item first,  smooth manifold limits diffeomorphic to
nearby manifolds (Anderson--Cheeger \cite{AC}, Petersen \cite{P});
\item next, topological manifold limits homeomorphic to nearby
manifolds (Grove--Peterson--Wu \cite{GPW});
\item then topologically stratified limits (Perel'man \cite{P1,P2},
Perel'man--Petrunin \cite{PP}); and finally
\item
more-singular limits currently not good for much.
\end{enumerate}
Limits in the homeomorphism case (2) were first only known to be
homology manifolds, and nearby  manifolds were analyzed using
controlled topology. However Perel'man \cite{P1} later used the
Alexandroff curvature structure to show the limits in (2) are in
fact manifolds, and extended this to stratifications with manifold
strata in some singular cases (3). This considerably simplified
the analysis and removed the need for homology manifolds. An
approach to singular cases using Ricci curvature is given by
Cheeger--Colding \cite{CC}, see also  Zhu \cite{Z}. We might still hope
for a role for the more sophisticated topology:
 \begin{problem} Are
there differential-geometric conditions or processes that give
exotic homology manifold limits?\end{problem}

Such conditions must involve something other than diameter, volume, and
sectional curvature bounds. Exotic ENR homology manifolds cannot arise
this way so the most interesting outcome would be to get
infinite-dimensional limits. An analysis of manifolds near such limits has
been announced by Dranishnikov and Ferry  and apparently these can vary
quite a lot, see also \fullref{S:nonENR}.
\subsection{Stratified Gromov--Hausdorff limits}
The most immediately promising problems about limits concern stratifications.

\begin{problem}\label{2.2} Are there differential-geometric conditions on
smooth stratified sets that imply Gromov--Hausdorff limits are homotopy stratified sets with homology-manifold strata? What is the structure of the nearby smooth stratified sets?\end{problem}

There are two phenomena here: ``collapse'' in which new strata are generated, and convergence that is in some sense stratum-wise.
The first case  is typified by ``volume collapse'' of manifolds to
stratified sets. The Cheeger--Fukaya--Gromov structure theory of collapsed manifolds suggests the nearby manifolds should be total spaces of stratified systems of fibrations with nilmanifold fibers. Note however this structure should only be topological: smooth structures on the limit and bundles are unlikely in general. In cases where curvature is bounded below Perel'man's analysis of the Alexandroff structure of the limit gives a topologically stratified space, solving the first part of the problem.

In the second case (stratum-wise convergence) the given smooth stratifications need not converge, but some sort of ``homotopy intrinsic'' stratifications should converge. More detail and an elaborate proposal for the topological part of this question is given in Quinn \cite{Q4}. In this case if strata in the limits are ENR homology manifolds then one expects nearby stratified sets to be stratified s-cobordant. Again cases where  Perel'man's  Alexandroff-space results apply should be much more accessible.

 The motivation for this question is to study compactifications of
 collections of algebraic varieties or of stratifications arising in
 singularity theory. Therefore to be useful the ``differential-geometric''
 hypotheses should have reasonable interpretations in these contexts.
 Other possibilities are to relate this to limits of special processes, eg the Ricci flow (Glickenstein \cite{G}) or special limits defined by logical constraints (van den Dries \cite{vD}).

\section{Compactifications}\label{S:compact}
Negatively curved spaces and groups (in the sense of Gromov) have
compactifications with ``boundaries'' defined by equivalence classes of
geodesics.  ``Hyperbolization'' procedures that mass-produce examples are
described by Davis--Januszkiewicz \cite{DJ},
Davis--Januszkiewicz--Weinberger \cite{DJW} and Charney--Davis \cite{CD}. ``Visible boundaries'' can be defined for nonpositively curved spaces using additional geometric information. Bestvina \cite{Ba} has given axioms for compactifications and shown compactifications of Poincar\'e duality groups  satisfying his axioms give homology manifolds.

 In classical cases the space on which the group acts is homeomorphic to Euclidean space, the boundary is a sphere compactifying the space to a disk. Behavior of limits in the sphere are more interesting  than the sphere itself. More interesting examples arise with ``Davis'' manifolds: contractible nonpositively curved manifolds not simply connected at infinity. Fischer \cite{Fr}  shows that a class of these have boundaries that are
\begin{enumerate}
\item finite dimensional cohomology manifolds with the (\v Cech) homology of a sphere;
\item not locally 1--connected (so not ENR);
\item homogeneous, and
\item the double of the compactification along the boundary is a genuine sphere.
\end{enumerate}
This connects nicely with the non-ENR questions raised in
\fullref{S:nonENR}. However it seems unlikely that interesting ENR
examples will arise this way:  boundaries can be exotic only if the input
space is exotic, for example the universal cover of an exotic closed
aspherical manifold as in \fullref{S:aspherical}, and this is probably not
compatible with nonpositive curvature assumptions, see \fullref{S:special}.

 To get more exotic behavior probably will require going outside the nonpositive curvature realm:

 \begin{problem}\label{3.1} Find non-curvature constructions for limits at infinity of  Poincar\'e duality groups, and find (or verify) criteria for these to be  homology manifolds.\end{problem}
      See Davis \cite{D} for a survey of Poincar\'e duality groups.  This question may provide an approach to  closed aspherical exotic homology manifolds: first construct the ``sphere at infinity'' of the universal cover, then somehow fill in.

 A variation on this idea is suggested by a proof of cases of the Novikov
 conjecture by  Farrell--Hsiang \cite{FH} and many others since. They use
 a compactification of the universal cover  to construct a fiberwise
 compactification of the  tangent bundle. This suggests directly
 constructing completions of a bundle rather than a single fiber. The
 bundle might include parameters, for instance  to resolve ambiguities
 arising in constructing limits without negative curvature. The context
 for this is discussed in \fullref{S:approx}.

   \begin{problem}\label{3.2} Construct ``approximate'' limits of duality groups, as ``fibers'' of the total space of an approximate fibration over a parameter space.\end{problem}

\section{Group actions and non--$\Z$ coefficients}\label{S:actions}
This topic probably has the greatest potential for profound applications,
but also has severe technical difficulty. Smith theory shows fixed sets of
actions of $p$--groups on homology manifolds must be $\Z/p\Z$ homology
manifolds. In the PL case there is a remarkable near converse:  Jones
\cite{Jl} shows  PL $\Z/p\Z$ homology submanifolds satisfying the Smith
conditions are frequently fixed sets of a $\Z/p\Z$ action.  Better results
are likely for topological  actions. \begin{problem}\label{4.1} Extend the
Jones analysis to determine when 1--LC embedded  $\Z/p\Z$ homology submanifolds are fixed sets of $\Z/p\Z$ actions.\end{problem}

If the submanifold is an ENR then there are tools available (eg mapping
cylinder neighborhoods) that should bring this within reach. Unfortunately
the non-ENR case is likely to be the one with powerful applications.  As a
test case we formulate a stable  version of \fullref{4.1} in which some difficulties should be avoided:

\begin{problem}\label{4.2} Suppose $X\subset R^n$ is an even-codimension properly embedded possibly non-ENR $\Z/p\Z$ homology manifold and is $\Z/p\Z$ acyclic. Is there a $\Z/p\Z$ action on $R^{n+2k}$ for some $k$ with fixed set $X\times\{0\}$?\end{problem}

Extending  \fullref{4.1} to a systematic classification theory for topological group actions will require a good understanding of  the corresponding homology manifolds:
\begin{problem}\label{4.3} Are there ``surgery theories'' for  $\Z/n\Z$ and  rational
homology manifolds? \end{problem}

Surgery for PL manifolds up to $\Z/p\Z$ homology equivalence was developed
in the 1970s by Quinn \cite{Q5}, Anderson \cite{A}, Taylor--Williams
\cite{TW}, and a speculative sketch for PL $\Z/p\Z$ homology manifolds
is given in Quinn \cite{Q5}. There are two serious difficulties for a
topological version.  The first problem is that the local and ``normal''
structures do not decouple. The boundary of a regular neighborhood in
Euclidean space is the appropriate model for the Spanier--Whitehead
dual of a space.  $\Z$--Poincar\'e spaces are characterized by this
neighborhood being equivalent to a spherical fibration over the
space. When manipulating a Poincar\'e space within its homotopy type
(eg  while constructing homology manifolds) the bundle gives easy and
controlled access to the Spanier--Whitehead dual. $\Z/p\Z$ Poincar\'e
spaces have regular neighborhoods that are $\Z/n\Z$ spherical fibrations,
but this only specifies the $\Z/p\Z$ homotopy type. Local structure at
other primes can vary from place to place, and the normal structure
must conform to this. Some additional structure is probably needed,
but this is unclear.

The second problem is that constructions of $\Z/p\Z$ homology manifolds
are unlikely to give ENRs. In the $\Z$  case ENRs are obtained as limits
of sequences of controlled homotopy equivalences. Homotopy equivalences
are obtained because obstructions to constructing these can be identified
with global data (essentially the topological structure on the normal
bundle). In the $\Z/p\Z$ case there will probably only be enough data to
get controlled $\Z/p\Z$  homology equivalences. It seems likely that
$n$--dimensional homology manifolds can be arranged to have covering
dimension $n$ and have nice $\bigl[\frac{n-1}2\bigr]$ skeleta, but  above the middle dimension infinitely generated homology prime to $p$ is likely to be common. 

$\Z/p\Z$ homology manifolds might geometrically implement some of the
remarkable but formal  ``$p$--complete'' manifold theory proposed in Sullivan \cite{S}.

Dranishnikov \cite{Dr2} gives constructions of rational homology 5--manifolds with large but still finite covering and cohomological dimension. If this complicates the development it may be appropriate to consider only homology manifolds with covering dimension equal to the duality dimension.

\subsection{ Circle actions}
A group action is ``semifree'' if points are either fixed by the whole
group or moved freely. In this case the fixed set is also the fixed set of
the $\Z/p\Z$ subgroups, all $p$, so it follows from Smith theory that it
is a $\Z$ homology manifold.  Problems \ref{4.1} and \ref{4.2} therefore
have analogs for semifree $S^1$ actions and $\Z$--coefficient homology manifolds:

  \begin{problem}\label{4.4} Determine when $\Z$ coefficient homology submanifolds  satisfying Smith conditions are fixed sets of semifree $S^1$ actions. \end{problem}
  This setting has the advantages that there are fewer obstructions, and
  in the ENR case the analog of \fullref{4.3} is already available.
Again the significance of non-ENR case depends on how many non-ENR
homology manifolds there are (see \fullref{S:nonENR}). If they all occur as boundaries with ENR interior then it seems likely a general action will be concordant to one with ENR fixed set. At the other extreme if there are fixed sets with non-integer local index then a  full treatment of group actions will probably need non-ENRs.

\section{ Non-ENR homology manifolds}\label{S:nonENR}
It is a folk theorem that a homology manifold that is finite dimensional
and locally 1--connected is an ENR. The proof goes as follows: duality
shows homology manifolds are homologically locally $n$--connected, all
$n$, and a local Whitehead theorem shows local 1--connected and homological
local $n$--connected implies local $n$--connected in the usual homotopic
sense. Finally finite-dimensional and locally $n$--connected for large $n$
implies ENR. The point here is that the ENR condition can fail in two
ways: failure of finite dimensionality or failure of local 1--connectedness. These are discussed separately.
 \subsection{ Infinite dimensional homology manifolds}
 Here we consider locally compact metric spaces that are locally
contractible (or at least locally 1--connected) and homology manifolds in the usual finite-dimensional sense, but with infinite covering dimension. This is not related to manifolds modeled on infinite dimensional spaces.
 
 \begin{problem}\label{5.1} Is there a ``surgery theory'' of infinite dimensional homology manifolds?\end{problem}

 These were shown to arise as cell-like images of manifolds by
 Dranishnikov \cite{Dr1}, following a proposal of Edwards. Recently
 Dranishnikov and Ferry have announced that there are examples with
 arbitrarily close (in the Gromov--Hausdorff sense) topological manifolds
 with different homotopy types. This contrasts with the ENR case where
 sufficiently close manifolds are all homeomorphic, and suggests this is a
 way to loosen the strait-jacket constraints of homotopy type in standard
 surgery. In particular the ``surgery theory'' should {\it not\/} follow
 the usual pattern of fixing a homotopy type, and ``structures'' should
 include manifolds of different homotopy type. This might be done by
 following Dranishnikov--Ferry in assuming existence of metrics that are
 sufficiently Gromov--Hausdorff close. See \fullref{3.1}.

 The source dimension for infinite-dimensionality is not quite settled:
 \begin{problem}\label{5.2} Are there infinite-dimensional $\Z$--homology
4--manifolds? Are there infinite-dimensional homology 4-- or 5--manifolds with nearby manifolds of different homotopy type?\end{problem}

Walsh \cite{W} has shown homology manifolds of homological dimension
$\leq3$ are finite dimensional. Dydak--Walsh \cite{DW} produced
infinite-dimensional examples of homology 5--manifolds but these do not
connect with the Dranishnikov--Ferry analysis. It may be that an interesting ``surgery theory'' does not start until dimension 6.

  \subsection{Non locally--1--connected homology manifolds}
Now consider finite dimensional metric homology manifolds that may fail to
be locally 1--connected. These arise as ``spheres at infinity'' for certain
groups, see \fullref{S:compact}. The first question seeks to locate these spaces relative
to ENR and ``virtual'' homology manifolds. This is important for
applications to group actions, see \fullref{4.4}. 

\begin{problem}\label{5.3} Extend the definition of local index to finite dimensional  non-ENR homology manifolds.\end{problem}
If the extended definition still takes values in $1+8\Z$ then the next
question (motivated by \fullref{S:compact}) would be:

\begin{problem}\label{5.4A} Does every finite dimensional homology manifold arise as the ``weakly tame'' boundary of one with ENR interior? Is the union of two such extensions along their boundary an ENR?\end{problem}

Here ``weakly tame'' should be as close to  ``locally 1--connected complement'' as possible. An affirmative answer to this question would suggest thinking of non-ENR homology manifolds as ``puffed up'' versions not much different from  ENRs.

At another extreme ``approximate'' homology manifolds are defined in section 6 in terms of approximate fibrations with homology manifold base and total space. These behave as though they have
``fibers'' that are homology manifolds with  local indices in  $1+8\Z_{(2)}$.  If an extension of the local index to non ENRs can take non-integer  values then the ENR boundary question above is wrong and we should ask:

\begin{problem}\label{5.4B} Does every finite dimensional homology manifold occur as a fiber of an approximate fibration with ENR homology manifold base and total space? Conversely is any such approximate fibration concordant to one with such a fiber?\end{problem}
An appropriate relative version of this would show approximate homology manifolds are equivalent to finite dimensional non-ENR homology manifolds.

\section{Approximate homology manifolds}\label{S:approx}
The intent is that approximate homology manifolds should be fibers of approximate fibrations with base and total space ENR homology manifolds. Actual point-inverses are not topologically well-defined and do not encode the interesting information, so we use a germ approach. For simplicity we restrict to the compact case (fibers of proper maps).

A {\it compact approximate homology manifold\/} is a pair $(f\co E\to B,b)$ where $f$ is a proper approximate fibration with homology manifold base and total space, and $b$ is a point in $B$.  ``Concordance'' is the equivalence relation generated by
\begin{enumerate}
\item changing the basepoint by an arc in the base;
\item restricting to a neighborhood of the basepoint; and
\item  product with identity maps of homology manifolds.
\end{enumerate}

\subsection{Basic structure} Suppose $F$ is a compact virtual homology manifold defined by an proper approximate fibration $f\co E\to B$  with $E, B$ connected homology manifolds, and $b\in B$.
\begin{enumerate}
\item $F$  has a well-defined homotopy type (the homotopy fiber of the map) that is a Poincar\'e space (with universal coefficients);
\item this Poincar\'e space has a canonical topological reduction of the normal fibration, given by restriction of the difference of the canonical reductions of $E$ and $B$; and
\item there is a local index  defined by $i(F)=i(E)/i(B)$.
\end{enumerate}
If $f$ is a locally trivial bundle then the fiber is a ENR homology
manifold. Multiplicativity of the local index shows the formula in (3)
does give the local index of the fiber in this case. In general the
quotients in (3) lie in $1+8\Z_{(2)}$, where $\Z_{(2)}$  is the rationals with odd denominator.

\subsection{ Example}
Suppose $X$ is a homology manifold, and choose a 1--LC embedding in a manifold of dimension at least 5. This has a mapping cylinder neighborhood; let $f\co \partial N\to X$ be the map. Duality shows this is an approximate fibration with fiber the homotopy type of a sphere. As an approximate homology manifold the index is $1/i(X)$, which is not an integer unless $i(X)=1$.

Products of these examples with genuine homology manifolds show that all
elements of  $1+8\Z_{(2)}$ are realized as indices of approximate homology manifolds.

\subsection{Example}
A tame end  of a manifold has an ``approximate collar'' in the sense of a neighborhood of the end that approximately fibers over $R$. In the controlled case the  local fundamental group is required to be stratified; see Hughes \cite{H} for a special case and Quinn \cite{Q6} in general.

There is a finiteness obstruction to finding a genuine collar. In some cases it follows that the approximate homology manifold appearing as the ``fiber'' of the approximate collar does not have the homotopy type of a finite complex.

\begin{problem}\label{6.1} Show that the exotic behavior of the examples are the only differences: if an approximate homology manifold has integral local index and is homotopy equivalent to a finite complex then it is concordant to an ENR homology manifold.\end{problem}

 \begin{problem}\label{6.2} Define ``approximate transversality'' of homology manifolds, and determine when a map from one homology manifold to another can be made approximate transverse to a submanifold. \end{problem}
The examples give maps approximately transverse to a point, but for which
more geometric forms of transversality are obstructed. Transversality
theories restricted to situations where indices must be integers have been
developed by Johnston  \cite{J}, Johnston--Ranicki  \cite{JR} and
Bryant--Mio  \cite{BM}, and a finiteness-obstruction case has been
investigated by Bryant--Kirby  \cite{BK}. The hope is that a more complete
approximate transversality theory is possible. There still will be
restrictions however: a degree--1 map of homology manifolds of different index cannot be made geometrically transverse to a point in any useful sense.

\begin{problem}\label{6.3} Develop a surgery theory for approximate homology manifolds.
\end{problem}
 The obstructions should lie in the $L^{-\infty}$ groups introduced by Yamasaki  \cite{Y}.

\section{Special metric spaces}\label{S:special}
Several special classes of metric spaces have been developed, particularly
by the Russian school, as general settings for some of the results of
differential geometry. It is natural to ask how these hypotheses relate to
manifolds and homology manifolds, but for the question to have much real
significance it is necessary to have sources of examples not a priori
known to be manifolds. Gromov--Hausdorff convergence gives Alexandroff
spaces, see \fullref{S:aspherical}.

 A Busemann space is a metric space in which geodesic (locally length-minimizing) segments can be extended, and small metric balls are cones parameterized by geodesics starting at the center point. The standard question is:
 
\begin{problem}\label{7.1} Must a Busemann space be a manifold?\end{problem}
This is true in dimensions $\leq4$; the 4--dimensional case was done by Thurston  \cite{T} and is not elementary.

An Alexandroff space is a metric space in which geodesics and curvature constraints make sense, but with less structure than  Busemann spaces. These need not be homology manifolds, so the appropriate question seems to be:

\begin{problem}\label{7.2} Is an  Alexandroff space that is a homology manifold in fact a manifold in the complement of a discrete set?\end{problem}
The problematic discrete set should be detectable by local fundamental
groups of complements, as with the cone on a non-simply-connected homology
sphere. The answer is ``yes''  when there is a lower curvature bound,
because the analysis in  Perel'man  \cite{P1,P2} and Perel'man--Petrunin  \cite{PP} shows it is a topological stratified set and topological stratified sets have this property (Quinn  \cite{Q7}).

\section{Low dimensions}\label{S:low}
1-- and 2--dimensional ENR homology manifolds are manifolds. In dimensions $\geq 5$ exotic homology manifolds of arbitrary local index exist, and there are ``many'' of them in
the senses that
\begin{itemize}
\item there is a ``full surgery theory'' given by
Bryant--Ferry--Mio--Weinberger  \cite{BFMW1} for dimensions $\geq6$ and
announced by Ferry--Johnston for dimension 5; and
\item in dimensions $\geq6$ Bryant--Ferry--Mio--Weinberger have announced a proof that an arbitrary homology manifold is the cell-like image of
one with the DDP.
\end{itemize}
The 5--dimensional case of (2) is still open:

\begin{problem}\label{8.1} Can 5--dimensional exotic homology manifolds be resolved by ones with the DDP?\end{problem}

In dimension 4:

\begin{problem}\label{8.2} Are there exotic 4--dimensional homology manifolds?\end{problem}

The expected answer is ``yes.'' In a little more detail the
possibilities seem to be: 
\begin{enumerate}
\item exotic homology manifolds
don't exist; or
\item sporadic examples exist; or
\item there is a ``full surgery theory''.
\end{enumerate}
Even in higher dimensions there are currently no methods for getting
isolated examples: to get anything one essentially has to go through the
full surgery theory. More-direct examples in higher dimensions would be
useful in approaching (2) as well as interesting in their own right. In
(3) note there is currently a fundamental group restriction in the
manifold case   Freedman--Quinn \cite{FQ}, Freedman--Teichner \cite{FT},
Krushkal--Quinn \cite{KQ}. Homology-manifold surgery would imply manifold surgery so ``full surgery theory'' should be interpreted to mean ``as full as the manifold case.''

Finally in dimension 3:

\begin{problem}\label{8.3} Are there exotic 3--dimensional homology manifolds?\end{problem}

The expected answer is ``no''. See Repov\v s \cite{R} for special conclusions in the resolvable case.

\section{ Homeomorphisms and the DDP}\label{S:DDP}
The basic question is: do the key homeomorphism theorems for
manifolds extend to homology manifolds? The question
should include a nondegeneracy condition that gives manifolds in
the index $=1$ case. Here we use the DDP (Disjoint Disk Property): any two
maps $i,j\co D^2\to X$  can be arbitrarily closely approximated by maps
with disjoint images. However see  \fullref{9.6}.

There is a feeling that the first three problems should be roughly equivalent in the sense that one good idea could resolve them all.

 \begin{problem}\label{9.1} Is the ``$\alpha$ approximation
theorem'' of
Chapman--Ferry  \cite{CF} true for DDP homology manifolds? \end{problem}

This expected answer is ``yes'', and current
techniques suggest a proof might break into two sub-problems:
\begin{itemize}\item A compact metric homology manifold $X$ has
$\epsilon>0$ so if $X'$ is $\epsilon$ homotopy equivalent to $X$
then there is a DDP $Y$ with cell-like maps onto both $X$ and
$X'$; and
\item if $X$, $X'$, or both, have the DDP then the corresponding
cell-like maps can be chosen to be homeomorphisms.
\end{itemize}
As a testbed for technique for the second part it would be
valuable to have a surgery-type proof of Edwards' theorem: a
cell-like map from a (genuine) manifold to a homology manifold
with DDP can be arbitrarily closely approximated by homeomorphisms.

 \begin{problem}\label{9.2} Is the h-cobordism
theorem true for DDP homology manifolds? \end{problem}

h-cobordisms appear in a natural way in the definition of ``homology manifold structure sets'', among other places,
and can be produced by surgery.

\begin{problem}\label{9.3} Is a homology
manifold with the DDP arc-homogeneous?\end{problem}

 ``Arc-homogeneous'' means if $x,y$ are in the same component of $M$ then there is a homeomorphism $M\times I\to M\times I$ that is the identity on one end and on the other takes $x$ to $y$. ``Isotopy-homogeneous" is the sharper version in which the homeomorphism is required to preserve the $I$ coordinate, so gives an ambient isotopy taking $x$ to $y$. The ``arc'' in the terminology refers to the track of the point under the homeomorphism or isotopy.

 An affirmative answer to 9.3 would show DDP homology manifolds have coordinate charts homeomorphic to subsets of standard models in the same way manifolds have Euclidean charts. Note this is consistent with a number of different models in each index:  only one model could occur in a {\it connected\/} homology manifold but different components might have different models.

  The ``Bing--Borsuk
conjecture'' is that a homogeneous ENR is a manifold, where ``homogeneous'' is used in the  traditional sense that any two points have
homeomorphic neighborhoods. A version more in line with current
expectations, and avoiding low dimensional problems, is that the
homogeneous ENRs of dimension at least 5 are exactly the DDP homology
manifolds. For applications and philosophical reasons we prefer arc
versions of homogeneity, and  split the question into  \fullref{9.3} and a homological question:

\begin{problem}\label{9.4} Is a locally 1--connected  homologically arc-homogeneous space a  homology manifold (possibly infinite-dimensional)?\end{problem}
A space is {\it homologically\/} arc-homogeneous if for any arc $f\co I\to X$ the induced maps
$$H_i(X\times\{0\},X\times\{0\}-(f(0),0);\Z)\to H_i(X\times I,X\times
I-\text{graph}(f);\Z)$$
is an isomorphism. This is clearly an analog of ``arc-homogeneous'' as
defined above. A homology manifold satisfies this by Alexander duality. In
fact it holds for $(I,0)$ replaced by a $n$--disk and a point in the boundary.

It was shown by Bredon
that homogeneous (in the traditional point sense) ENRs are homology manifolds {\it provided\/} the local homology groups are finitely generated, see Bryant  \cite{B1},
Dydak--Walsh  \cite{DW}.  The problem is to show the local homology groups
form a locally constant sheaf. The arc version of homogeneity gives local
isomorphisms so the problem becomes showing these are locally
well-defined. This would follow immediately from a ``homologically
2--disk-homogeneous'' hypothesis, so is equivalent to this condition. The
question is whether this follows from arc-homogeneity and local
1--connectedness. Bryant  \cite{B2} has recently proved \fullref{9.4}
under the assumption that  the space is an ENR. Note that a finite
dimensional locally 1--connected homology manifold is an ENR, so the question remaining is whether ``ENR'' can be shifted from hypothesis to conclusion.

In a somewhat different direction the following is still unknown even in the manifold case:

\begin{problem}\label{9.5} Is the product of a homology manifold and $R$ homogeneous? \end{problem}

This can be disengaged from the homogeneity questions by asking
``does $X\times R$ have DDP?'', but see the discussion of the DDP
in \fullref{9.6}.
$X\times R^2$ does have DDP (Daverman \cite{D2} and there are quite a
number of properties of $X$  that imply
$X\times R$ has DDP, see
Halverson \cite{H}, Daverman--Halverson \cite{DH}. However there are ghastly (in the
technical sense) examples of homology manifolds that show none of
these  properties holds in general, see Daverman--Walsh \cite{DW},
Halverson \cite{H}.

The final question is vague but potentially important:

\begin{problem}\label{9.6} Is
there a weaker condition than DDP that implies index $=1$ homology
manifolds are manifolds?\end{problem}

 If so then this condition should be
substituted for the DDP in the other problems in this section.
This could make some of them significantly easier, and may also
help with understanding dimension 4. A good way to approach this
would be to find a surgery-based proof of Edwards' approximation
theorem (see 9.1), then inspect it closely to find the minimum
needed to make it work. Edwards' proof (see Daverman \cite{D1}) uses
unobstructed cases of engulfing and approximation theorems.
Surgery by contrast  proceeds by showing an obstruction vanishes.
Potentially-obstructed proofs (when they work) are often more
flexible and have led to sharper results.

\bibliographystyle{gtart}
\bibliography{link}

\end{document}